\documentclass{article}
\pdfoutput=1
\usepackage{amssymb,amsfonts,amsmath,amsthm}
\usepackage{epsfig}
\usepackage{longtable,verbatim}

\newtheorem{theorem}{Theorem}[section]{\bfseries}{\itshape}
\newtheorem{Prop}[theorem]{Proposition}{\bfseries}{\itshape}
\newtheorem{Lem}[theorem]{Lemma}{\bfseries}{\itshape}
\newtheorem{Cor}[theorem]{Corollary}{\bfseries}{\itshape}
{\bfseries}{\itshape}
{\bfseries}{\itshape}
\newtheorem{definition}[theorem]{Definition}{\bfseries}{\itshape}
\newtheorem{remark}[theorem]{Remark}{\bfseries}{\itshape}

\long\def\symbolfootnote[#1]#2{\begingroup\def\thefootnote{\fnsymbol{footnote}}
\footnote[#1]{#2}\endgroup}

\newcommand{\sgn}{\ensuremath{\mathrm{sgn}}}
\newcommand{\finfldq}[1]{\ensuremath{\mathbb{F}_{#1}}}
\newcommand{\polyringq}[1]{\ensuremath{\mathbb{F}_{#1}[t]}}

\newcommand{\ints}{\ensuremath{\mathbb{Z}}}

\newcommand{\ratsfuncq}[1]{\ensuremath{\mathbb{F}_{#1}(t)}}
\newcommand{\disc}[1]{\ensuremath{\mathrm{disc}(#1)}}
\newcommand{\ratsfuncqextn}[2]{\ensuremath{\mathbb{F}_{#1}(t,#2)}}
\newcommand{\charfld}{\ensuremath{\mathrm{char}}}
\newcommand{\norm}[1]{\ensuremath{\mathcal{N}(#1)}}
\newcommand{\algints}{\ensuremath{\mathcal{O}_K}}

\title{Simple cubic function fields and class number computations}

\author{Pieter Rozenhart \\ Department of Mathematics, University of Calgary\\
 2500 University Drive NW\\
  Calgary, Alberta, T2N 0J1 Canada \\ \texttt{pmrozenh@alumni.uwaterloo.ca} \and Jonathan Webster \ \\   Department of Mathematics, Bates College, \\ 3 Andrews Road,
  Lewiston, ME 04240, USA\\
   \texttt{jwebster@bates.edu}}

\begin{document}
\maketitle

\begin{abstract}
In this paper, we study simple cubic fields in the function field setting, and also generalize the notion of a set of exceptional units to cubic function fields, namely the notion of $k$-exceptional units.  We give a simple proof that the Galois simple cubic function fields are the immediate analog of Shanks simplest cubic number fields.  In addition to computing the invariants, including a formula for the regulator, we compute the class numbers of the Galois simple cubic function fields over $\finfldq{5}$ and $\finfldq{7}$ using truncated Euler products.  Finally, as an additional application, we determine all Galois simple cubic function fields with class number one, subject to a mild restriction.
\end{abstract}

\symbolfootnote[0]{Keywords:  cubic function fields, ideal class number, regulator, computational number theory}

\section{Introduction}

By studying cubic number fields with exceptional units, computational number theorists can consider a class of number fields whose regulators are easy to calculate and this leads to a simple computation of the class number.   The quintessential example of this was Daniel Shanks' cubic number fields \cite{ShanksSimp} arising from the minimal polynomial $f(x)= x^3 - ax^2 - (a + 3)x - 1$.  These extensions are Galois and he imposed the further restriction that $a^2 + 3a + 9$ be prime.     With the regulator in hand and knowledge of the splitting behavior of all finite places, we may compute the class number in a straightforward manner.  He called these number fields the simplest cubic number fields.  The quantity $a^2 + 3a + 9$ was chosen to be prime so that the splitting of the finite places would be easy to calculate.  If we have a prime factorization of $a^2 + 3a + 9$, then the splitting of the finite places is again easy to calculate and computing the class number proceeds in the same computational manner as before.  

Two fundamental units in a cubic number field are said to be \emph{exceptional} if they differ by $1$.  Shanks' family gives rise to Galois cubic number fields with exceptional units.  Nagell \cite{Nagell1} showed that there are only two classes of cubic number fields with this property.  Ennola \cite{Ennola} studied the behavior of the class number for the noncyclic family that has minimal polynomal $g(x) = x^3 + kx^2 - (k + 1)x - 1$.  Much like the simplest cubic number fields, these number fields have regulators that are easy to calculate and so the class number may be obtained with relatively little work.  We refer the reader to \cite{Ennola} for a better list of references for the number field case.  

In this paper, our goal is to study the function field (over a finite field \finfldq{q}) analogue of these polynomials.  While several authors have considered the immediate analogue of Shanks' polynomial, including \cite{GaalPohst,PacRose}, the lack of a detailed study of these cubic function fields is the primary motivation for this work.    In the function field setting, one might ask: ``are there other starting points?"  The added freedom of having $q-1$ units might lead one to believe that by considering function fields that there could be more classes of ``simplest cubic (function) fields".   Our answer is mixed.  In the case of the direct analog of the Shanks' family, we answer no.  We prove that if a function field has the properties of a Galois simple cubic (function) field then it has the Shanks polynomial as a defining polynomial.  The proof is similar to Nagell's \cite{Nagell1} in the number field case.  Our proof recovers  Kersten and Michali\v{c}ek's \cite{Haberle,KerstMicha} results, but our method is much simpler and more elementary, while also taking into account the units of ring of algebraic integers.  However, by considering a wider class of simple cubic function fields that includes non-Galois fields, we show that there are many more than two families.  

We also prove a number of properties of the Shanks simple cubic function fields, including computation of the invariants of such fields.  Unlike purely cubic function fields where the regulator tends to be very large \cite{LeeSchYar} (exponentially large in the size of the field), the each member of the Shanks family has regulator that is both easy to compute and polynomial in the size of the parameter $A$ given in the Shanks polynomial.  In addition, we compute the class numbers of these fields up to a certain genus level over $\finfldq{5}$ and $\finfldq{7}$, using approximations of $L$-polynomials as described in \cite{SchStein1}.  As an application of our computation and the Hasse-Weil bound, we determine all simple cubic function fields with ideal class number equal to one with $A^2 + 3A +9$ cube-free.  Similar results for Shanks simple cubic number fields were obtained by Lettl \cite{Lettl} for class number one and Byeon \cite{Byeon2} for class number three simple cubic number fields, respectively.

The paper is organized as follows.  We review some preliminaries from the theory of algebraic function fields in 
Section \ref{sec:prelim}.  The main theoretical results, including the computation of invariants, appears in Section \ref{sec:FieldProps}.  Special attention is paid to proving that our units
are fundamental in Section \ref{sec:Units}, where we also compute a simple formula for the regulator.
  We give a brief summary of \cite{SchStein1} and the results of our computations in Section \ref{sec:compres}.  As an application of our computational results, we determine all simple cubic function fields with class number one, with $A^2 + 3A + 9$ cube-free.  We conclude with some open problems and future research directions in Section \ref{sec:concl}.

\section{Preliminaries}\label{sec:prelim}

For a complete introduction to algebraic function fields, see Rosen \cite{Rosen} or Stichtenoth \cite{Sticht}.    Let
\finfldq{q}\ be a finite field of characteristic at least~$5$, and
set $\finfldq{q}^* = \finfldq{q} \backslash \{0\}$. Denote by
\polyringq{q}\ and \ratsfuncq{q}\ the ring of polynomials and the
field of rational functions in the variable $t$ over \finfldq{q},
respectively.    For any non-zero $P \in \polyringq{q}$, let $\deg(P)$ denote its degree, and denote by $\sgn(P)$ the leading coefficient of $P$. For $P = 0$, we set $\sgn(P) = 0$. These notions extend in the obvious way to \ratsfuncq{q}. 

An \emph{algebraic function field} is a finite extension $K$ of $\ratsfuncq{q}$; its degree is the
field extension degree $n = [K : \ratsfuncq{q}]$. It is always possible to write a function
field as $K = \ratsfuncqextn{q}{y}$ where $F(t, y) = 0$ and $F(Y )$ is a monic polynomial of degree
$n$ in $Y$ with coefficients in $\polyringq{q}$ that is irreducible over $\ratsfuncq{q}$.   We assume that $\finfldq{q}$ is the full constant field of $K$, i.e. $F(Y)$ is absolutely irreducible.

The powers $y^i$ , $0 \leq i \leq n -1$, form an $\ratsfuncq{q}$-basis of the $\ratsfuncq{q}$-vector space $K$.  The \emph{$n$ conjugate mappings} map $y$ to the $n$ (distinct)
roots $y = y^{(0)}, y^{(1)}, . . . , y^{(n-1)}$ of $F(Y)$. Extending these mappings $\ratsfuncq{q}$-linearly to
$K$ now defines for every $\alpha \in K$ its $n$ conjugates $\alpha = \alpha^{(0)}, \alpha^{(1)}, \ldots , \alpha^{(n-1)}$.  Since we only consider cubic fields in this paper, we will frequently use the notation $\alpha,\alpha',\alpha''$ for the three conjugate roots of an irreducible cubic polynomial.

The discriminant of $n$ elements $\alpha_0, \alpha_1, . . . \alpha_{n-1} \in K$ is
\[ 	\disc{ \alpha_0,\ldots,\alpha_{n-1}} = \det (\alpha_i^{(j)})_{0 \leq i,j \leq n}^2 \in \ratsfuncq{q}.\]
If $\alpha_i = \alpha^i$ for some non-zero $\alpha \in K$ and $0 \leq i \leq n-1$, then $\disc{\alpha} =
\disc{1, \alpha, . . . , \alpha^{n-1}}$ is simply called the discriminant of $\alpha$, and
\[\disc{\alpha} = \prod_{i<j}(\alpha^{(i)} - \alpha^{(j)})^2.\]
We have $\disc{y} = \disc{F}$, the discriminant of $F$ as a polynomial in $Y$.

The \emph{maximal order} or \emph{coordinate ring} \algints\ of $K/\ratsfuncq{q}$ is the integral closure of $\polyringq{q}$ in $K$, and is analogous to the ring of algebraic integers of an algebraic number field. It is a free $\polyringq{q}$-module of rank $n$, and an $\polyringq{q}$-basis of \algints\ is an \emph{integral basis} of $K/\ratsfuncq{q}$. The discriminant of $K/\ratsfuncq{q}$
is $\disc{K} = \disc{\alpha_0, \alpha_1, . . . , \alpha_{n-1}}$ where $\{\alpha_0, \alpha_1, . . . , \alpha_{n-1}\}$ is any integral basis of $K/\ratsfuncq{q}$.  The polynomial $\disc{K} \in \polyringq{q}$ is independent of the basis chosen and unique up to square factors in $\finfldq{q}^*$ . For every non-zero element $\alpha \in K$, the \emph{index} of $\alpha$, denoted by $I(\alpha)$ or $I$ and unique up to square factors $\in \finfldq{q}^*$, is the rational function in $\ratsfuncq{q}$ satisfying $\disc{\alpha} = I^2 \cdot \disc{K}$. If $\alpha \in \algints$, then $I \in \polyringq{q}$, so $\disc{\alpha} \in \polyringq{q}$.

We note that $I \in \finfldq{q}^*$ , i.e., $\disc{y} = \disc{K}$ up to constant square
factors, if and only if the powers $y^i$ , $0 \leq i \leq n-1$, form an integral basis of $K/\ratsfuncq{q}$.
That is, if and only if $\algints = \finfldq{q}[t][y] = \finfldq{q}[t,Y]/(F(Y))$ where $(F(Y))$ is the principal
ideal generated by $F(Y)$ in $\finfldq{q}[t,Y]$.

A \emph{place} of a function field $K$ is simply the maximal ideal of some valuation ring of $K$.
The places of $\ratsfuncq{q}$ consist of the finite places, identified with the monic irreducible
polynomials in $\polyringq{q}$, and the place at infinity $P_\infty$, identified with the rational function
$1/x$.  For any finite place $P$ of \ratsfuncq{q}, the degree $\deg(P)$ is the degree of the polynomial $P$, and the degree of the infinite place $\deg(P_\infty)$ is set to be 1. The completion of $\ratsfuncq{q}$ with respect to any
place $P$ is the field of \emph{Laurent series} $\finfldq{q^d} \left\langle P \right\rangle$ in $P$ over \finfldq{q^d}, where $d = \deg(P)$. Non-zero elements in this field have the form $\sum_{i \geq m}a_iP^i$ where $m \in \ints, a_i \in \finfldq{q^d}$ for $i \geq m$, and $a_m \neq 0$.  The notions of degree and sign can be extended to Laurent series, namely we write $\deg(\alpha) = -m$ and $\sgn(\alpha) = a_m$ for all $\alpha \in \finfldq{q^d}\left \langle P \right \rangle$.

For $P = P_\infty$, we have $\finfldq{q^d}\left \langle P \right\rangle = \finfldq{q}\left \langle t^{-1} \right \rangle$.    Note that $ \finfldq{q}\left \langle t^{-1} \right \rangle$ is the completion with respect to the infinite place of $\ratsfuncq{q}$.  We can view the roots of any cubic polynomial having coefficients in $\polyringq{q}$ as lying in this field of Laurent series.  We make use of this fact briefly in Section~\ref{sec:Units}. 

If $K$ is any finite algebraic extension of $\ratsfuncq{q}$ of degree $n$, then the place at infinity of $\ratsfuncq{q}$ splits in $K$ as $P_\infty = \mathfrak{p_1}^{e_1}\mathfrak{p_2}^{e_2} \cdots \mathfrak{p_s}^{e_s}$, where $s \in \mathbb{N}$ and for $1 \leq i \leq s$, $\mathfrak{p_i}$ is a place of $K$ of residue degree $f_i \in \mathbb{N}$. and ramification index $e_i \in \mathbb{N}$ with $\sum_{i=1}^{s} e_i f_i = n$.  If we sort the pairs $(e_i, f_i)$, $1 \leq i \leq s$, in lexicographical order, then the $2s$-tuple $(e_1, f_1, e_2, f_2, . . . , e_s, f_s)$ is the \emph{signature} of $K/\ratsfuncq{q}$.  If one replaces the place at infinity $P_\infty$ with any finite place, the tuple is referred to as the \emph{$P$-signature}.

Associated to a function field $K$ are its \textit{divisor class number} $h$, which is simply the order of the group $\mathcal{J}$ of $\finfldq{q}$-rational points on the Jacobian of $K$.  The ideal class group $\mbox{Cl}(\algints)$ is the factor group of fractional ideals modulo principal ideals in \algints. Its order is the \textit{ideal class number} $h'$ of $K$.

The \emph{unit group of $K/F$} is the group of units $\algints^*$ of the ring of integers $\algints$ of $K/\ratsfuncq{q}$.  As per Dirichlet's Unit Theorem, the group $\algints^*$ is an infinite Abelian group whose torsion part is $\finfldq{q}$ and whose free part has rank $s-1$, where $s$ is the number of places at infinity in $K/\ratsfuncq{q}$.  A set of generators for the free part of the unit group is called a \emph{system of fundamental units}.  The cubic function fields under consideration in this paper all have unit rank two, so we specialize to this case immediately.  If $\varepsilon_1$ and $\varepsilon_2$ are two independent fundamental units then the \emph{regulator} $R$ of the function field is the absolute value of the determinant of \[  \left( \begin{array}{cc}   \deg(\varepsilon_1^{(i)}) &  \deg(\varepsilon_1^{(j)}) \\  \deg(\varepsilon_2^{(i)}) &  \deg(\varepsilon_2^{(j)}) \end{array} \right)  \]
where $i,j \in \{0,1,2\}$ and $i \neq j$.

When the unit rank is greater than 0, the equality (see \cite{Sch1}, Section 5)
\[ h = \frac{R}{f}h'\]
holds, where $f$ is the common divisor of the degrees of all the infinite places and $R$ is the regulator.  The fields considered in this paper have signature $(1,1;1,1;1,1)$, so in these cases, it follows that $h = Rh'$.

We consider an absolutely irreducible nonsingular affine plane curve $\mathcal{C}_0$ defined by an equation $H(t,x)= 0$ where $H \in \mathbb{F}_q[t][x]$ is a bivariate polynomial of degree 3 in $x$.   Since we want to consider cubic function fields where one of the roots of its minimal polynomial is a fundamental unit, we may restrict our consideration to bivariate polynomials of the following form:
\[ H(t,x) = x^3 + Px^2 + Ax + c \]
where $P, A \in \mathbb{F}_q[t]$ and $c \in \mathbb{F}_q^*$.  We will consider the function field $K = \mathbb{F}_q(t,y)$ where $H(t,y) = 0$.  With these restrictions $K$ is a cubic function field with a root of the polynomial $H$ being a non-torsion unit in $\algints$.  We also require that the other fundamental unit differ from this root by a unit in $\finfldq{q}^*$.  We will prove that these two non-torsion units
are fundamental units later in the paper. 

\begin{definition}
Let $K$ be a cubic function field of unit rank two.  A pair of units $\epsilon$ and $\tau$ in $K$ are said to be \emph{$k$-exceptional} if $\tau = \epsilon + k$ for some $k \in \finfldq{q}^*$.  A pair of 1-exceptional units will simply be called a pair of exceptional units.
\end{definition}

Hereafter, a cubic function field having $k$-exceptional units will be called a \emph{simple cubic function field}.  We will sometimes refer to a Galois simple cubic function field as a \emph{Shanks simple cubic function field}.
A cubic function field having $k$-exceptional units places restrictions on the $P$ and $A$ allowed to define $H(x,t)$.  Suppose that $\varepsilon$ is a unit of a cubic function field $K$ of unit rank two, with $\varepsilon \notin \finfldq{q}^*$.  Suppose further that $\varepsilon$ is a root of the polynomial 
\[H(t,x) =  x^3 + Px^2 + Ax + c,\]
where $c \in \finfldq{q}^*$, and $P,A \in \polyringq{q}$ and $\deg(A), \deg(P) \geq 1$, in order to avoid constant field extensions.  
Suppose that $\varepsilon + k$ is also a fundamental unit in $K$, with $\varepsilon + k \notin \finfldq{q}^*$ and $k \in \finfldq{q}^*$.  Then the norm of $\varepsilon + k$ is 
\begin{eqnarray*}
\norm{\varepsilon + k} & = & (\varepsilon + k)(\varepsilon + k)'(\varepsilon + k)'' \\
                                   & = & (\varepsilon' + k)(\varepsilon' + k)(\varepsilon'' + k) \\
                                   & = & \varepsilon \varepsilon' \varepsilon'' + (\varepsilon + \varepsilon'  + \varepsilon'') + (\varepsilon\varepsilon' + \varepsilon \varepsilon'' + \varepsilon' \varepsilon'') + k^3 \\
                                   & = & \norm{\varepsilon} + (\varepsilon + \varepsilon'  + \varepsilon'') + (\varepsilon\varepsilon' + \varepsilon \varepsilon'' + \varepsilon' \varepsilon'') + k^3\\
                                   & = & A - P - c + k^3.
\end{eqnarray*}

Since $\norm{\varepsilon + k} \in \mathbb{F}_q^*$, let $ n = \norm{\varepsilon+k}$.  Now a simple cubic function field has the following parameterization
\[ H(t,x) = x^3 + Px^2 + (P + c + n - k^3)x + c. \]
This may be thought of as a three parameter family involving the quantities $P$, $c$, and $c + n - k^3$.  Note this differs from the number field case, where for a given $P$ one could construct at most two simple cubic number fields (one cyclic and one not).  For instance, letting $P = t^3 \in \polyringq{5}$, we see that there are 16 different families of simple cubic function fields to consider.  

\section{Simple Cubic Function Fields}\label{sec:FieldProps}

In this section, we prove that cyclic cubic extensions of $\ratsfuncq{q}$ with $k$-exceptional units are the Shanks' family \cite{ShanksSimp}.  In particular, if $\varepsilon$ and $\varepsilon + k$ are fundamental units of a cyclic cubic function field with $N(\varepsilon) = c$ and $N(\varepsilon + k) = n$ then $c = -1$, and $k^3 - n = 2$ where $n,c,k \in \finfldq{q}^*$.  Part of the proof is similar to Nagell's \cite{Nagell1} for cubic number fields.  It is interesting that the additional units in $\finfldq{q}^*$
do not provide any new families of cyclic simple cubic function fields.  

\begin{theorem}
The cyclic cubic function fields of characteristic greater than five of unit rank two with $k$-exceptional units are given by the minimal polynomial $x^3 - Ax^2 - (A+3)x -1$, with $A \in \polyringq{q}$ and $\deg(A) \geq 1$, which is the function field analogue of Shanks' family of the simplest cubic number fields.
\label{simpleFieldConstruction}
\end{theorem}

\begin{proof}
Suppose that $\varepsilon$ is a unit of a cubic function field $K$ of unit rank two, with $\varepsilon \notin \finfldq{q}^*$.  Suppose further that $\varepsilon$ is a root of the polynomial and let
\[H(t,x) =  x^3 + Px^2 + (P + c + n - k^3)x + c\]
be as before.  The discriminant of $\varepsilon$ is simply the discriminant of $H(t,x)$, namely
\begin{eqnarray*}
D(\varepsilon)  & = & P^4 + (2n -2c - 2k^3 -4)P^3 \\
                     & & + ( 6c - 12n + 12k^3 + (k^3-c-n)^2 )P^2\\
                     & & + ( 18c( c+ n -k^3) - 12(k^3-c-n)^2)P\\
                     & & 4(k^3 - c - n)^3 + 27c^2 
\end{eqnarray*}
Thus, $D(\varepsilon)$ is a square if and only if $D(\varepsilon) = (P^2 + aP + b)^2$ with $a,b \in \finfldq{q}$, where $2a = 2n -2c - 2k^3 -4$, $2b+a^2 = 6c - 12n + 12k^3 + (k^3-c-n)^2$, $2ab = 18c( c+ n -k^3) - 12(k^3-c-n)^2$ and $b^2 = 4(k^3 - c - n)^3 + 27c^2$.  

We note first that the parameters  $c = -1$, and $k^3 - n = 2$ always give rise to the minimal polynomial $x^3 - Px^2 + (P-3)x -1$ with discriminant $(P^2 -3P + 9)^2$.  That is, since $2a = 2n - 2c - 2k^3 -4$, then $a = n - k^3 - c - 2 = -3$.  Since $2b+a^2 = 6c - 12n + 12k^3 + (k^3-c-n)^2$, we can write $b = (c - 2)(1 -2k^3 + 2n) = 9$.  We now show that there are no other solutions.  In doing so, we use the identities $2ab = 18c( c+ n -k^3) - 12(k^3-c-n)^2$ and $b^2 = 4(k^3 - c - n)^3 + 27c^2$.

Using the identities obtained for $a$ and $b$ to write $2ab = 18c( c+ n -k^3) - 12(k^3-c-n)^2$ completely in terms of $c, k, n$ we obtain
\[ 0 = ( c + 1 ) (n - k^3 + 2 ) ( (n - k^3 + 2) - (c + 1) ) = CB(C-B) .\]
with $C = n - k^3 + 2 $ and $B = c + 1 $.

Using the $k,c,n$ as in the theorem statement gives $C = B = C-B = 0$.  Naturally, if any two of $C, B, C-B$ are chosen to be zero, then the third is as well.  As seen above, if another family of fields with $k$-exceptional units are to be found, it will occur only by allowing exactly one of $C, B, C-B$ to be zero.  There are three cases to be considered and each is treated similarly.  We choose exactly one of $\{C,B,C-B\}$ to be zero and use $b^2 = 4(k^3 - c - n)^3 + 27c^2$ to derive a contradiction.  

Assume $ c + 1 \neq 0$ and $n - k^3 + 2 = 0$.  Substituting $k^3 - n = 2$ and $ b= (c - 2)(1 -2k^3 + 2n)$ into $b^2 = 4(k^3 - c - n)^3 + 27c^2$, gives a cubic expression in $c$.  We have $ 0 = -4 ( 1 + 3c + 3c^2 + c^3) = -4( 1 + c)^3$.  However, we assumed $c + 1 \neq 0$, yielding a contradiction.

Assume $c + 1 = 0$ and $n - k^3 + 2 \neq 0$.  Substituting $c = -1$ and $b= (c - 2)(1 -2k^3 + 2n)$ into $b^2 = 4(k^3 - c - n)^3 + 27c^2$, gives a cubic expression in $k^3 -n$.  We have $ 0 = 4( (k^3 -n)^3 -6(k^3-n)^2 + 12(k^3-n) -8) = 4( k^3 -n - 2)^3$ which once again gives the desired contradiction.

Finally, assume $c + 1 = n - k^3 + 2$ and $c + 1 \neq 0$.  Substituting $k^3 -n = 1 - c$ and $b= (c - 2)(1 -2k^3 + 2n)$ into $b^2 = 4(k^3 - c - n)^3 + 27c^2$, gives a quartic expression in $c$.  We have $0 = -4c(c+1)^3$.  Since $c$ has to be nonzero and $c + 1 \neq 0$ this gives the desired contradiction.

Hence, our family of polynomials is given by $f(x)= x^3 - Px^2 - (P-3)x - 1$.  Replacing $P$ with $A+3$, and considering the polynomial $x^3f(-1/x)$ yields the same field, defined by $x^3 - Ax^2 - (A+3)x -1$.  This completes the proof of the theorem.
\end{proof}

 From the proof of Theorem \ref{simpleFieldConstruction}, it also follows that many different polynomial families can be obtained in the non-Galois case.  This does not occur in the number field case.  
 We state and prove a number of properties of Shanks' family of cubic function fields.


\begin{Prop}
Let $f(x) = x^3 - Ax^2 -  (A+3)x - 1$, with $A \in \polyringq{q}$, $A^2 + 3A + 9$ cube-free and $K = \ratsfuncqextn{q}{\varepsilon}$, where $\varepsilon$ is a root of $f(x)$.  Then the following properties hold:
\begin{enumerate}
\item The discriminant $D(K)$ of $K$ is $D(K) = (A^2 + 3A + 9)^2/I^2$ where $I = \gcd(A^2+3A+9,A')$, and an integral basis of $\algints$ is $\{1,\theta,\theta^2/I\}$, where $\theta = \varepsilon-\frac{1}{3}A$.
\item The genus $g$ of $K$ is given by $g = 2\deg(A) - \deg(I) - 2$.
\item The field $K$ has signature $(1,1;1,1;1,1)$.
\end{enumerate}
\label{ListOfInvar}
\end{Prop}

\begin{proof}
By replacing $x$ with $x - \frac{1}{3}A$ in $f(x)$, we eliminate the quadratic term and obtain
\begin{equation}g(x) = x^3 - \frac{1}{3}(A^2+3A + 9)x - \frac{1}{27}(2A + 3)(A^2 +3A + 9). \label{stdform}\end{equation}
Set $\mathcal{A} =  \frac{1}{3}(A^2+3A + 9)$ and $\mathcal{B} = - \frac{1}{27}(2A + 3)(A^2 +3A + 9)$.  Then $\theta$ is a root of $g$.  The assumption that $A^2 + 3A + 9$ is cube-free is needed so that equation (\ref{stdform}) is in standard form (i.e.\ there is no non-constant polynomial $Q \in \polyringq{q}$ with $Q^2 \mid \mathcal{A}$ and $Q^3 \mid \mathcal{B}$). 

Since $(A^2 + 3A + 9)$ divides both $\mathcal{A}$ and $\mathcal{B}$, and $2A+3$ never divides $A^2 +3A + 9$ when $\charfld(\finfldq{q}) \neq 3$, it follows that $2A+3$ and $A^2 +3A + 9$ are relatively prime.  It is now simply a matter of checking if $(A^2 + 3A + 9)$ is square-free to find $I$.  This is done by computing $\gcd(A^2+3A+9,(A^2 + 3A + 9)') = \gcd(A^2 + 3A + 9,A'(2A+3)) = \gcd(A^2 +3A + 9,A')$, with the last equality following from $\gcd(A^2 +3A + 9,2A+3)=1$.  The result on the integral basis follows immediately upon applying Theorem 6.4 of \cite{LRWWS} with $U, T$ and $V$ all zero.

For part (2) or the Proposition, we use the Hurwitz genus formula. By the Hurwitz genus formula \cite{Sticht}, it follows that $g = \frac{1}{2}(\deg(D(K)) + \epsilon_{P_\infty}(K))-2$, where
\begin{equation}
\epsilon_{P_\infty}(K) =
 \left\{ \begin{array}{ll}
                                                2 & \mbox{if $3 \deg(\mathcal{A}) < 2 \deg(\mathcal{B})$ and $\deg(\mathcal{B}) \equiv 0 \pmod{3}$,} \\ 
                                                1 & \mbox{if $\deg(D(K))$ is odd,} \\ 
                                                0 & \mbox{otherwise,} 
                                            \end{array} \right.
\end{equation}
as in equation (3.2) of Landquist et al.\ \cite{LRWWS}.  Note that $\deg(\mathcal{B}) = 3 \deg(\mathcal{A})$, and that $\deg(D(K))$ is always even.  Hence, $\epsilon_{P_\infty}(K) = 0$, and since $\deg(D(K)) = 4 \deg(A) - 2 \deg(I)$, the result follows.

For the third property, note that $3 \deg(\mathcal{A}) = 3(\deg(A^2 + 3A + 9)) = 6 \deg(A)$, and $2 \deg(\mathcal{B}) = 2(3\deg(A)) = 6 \deg(A)$.  Furthermore, $\deg(D(K))$ is a square and so $\sgn(D(K))$ is a square in $\finfldq{q}$.
Hence, via the signature characterization of cubic fields given in \cite{LRWWS}, it follows that $K$ either has signature $(1,3)$ or $(1,1;1,1;1,1)$.  

Using the equivalent signature characterization in \cite{Sch1}, we note that $4\,\sgn(\mathcal{A})^3$ $= (4/27)(\sgn(A))^6 = 27\, \sgn(\mathcal{B})^2$, finally yielding $(1,1;1,1;1,1)$ as the signature of the field $K$.
\end{proof}

\begin{remark}
 The assumption that $A^2 + 3A + 9$ is cube-free is also helpful because it prevents the consideration of certain isomorphic function fields.  The fact that the fields that arising from  $f(x) = x^3 - Ax^2 -  (A+3)x - 1$ and $h(x) = x^3 - A^px^2 -  (A+3)^px - 1$ are isomorphic follows by properties of the Frobenius map. 
\end{remark}

In practice, $A^2 + 3A + 9$ is expected to be cube-free.  If the polynomial viewed as a random polynomial, then we expect it to be cube free with probability $1-1/q^3$.  We note that $\deg(I) \leq \deg(A') \leq \deg(A) - 1$ when $A^2 - 3A - 9$ is cube-free.  We can use this in concert with the genus formula to bound the size of $\deg(A)$ in terms of $g$.  This allows us consider finitely many $A$ for a fixed genus level in order to determine all simple cubic function fields with $A^2 + 3A + 9$ cube-free and class number $1$ in Section \ref{sec:compres}.  This is contained in the following lemma.

\begin{Lem}
Let $f(x) = x^3 - Ax^2 -  (A+3)x - 1$, with $A \in \polyringq{q}, A \notin \finfldq{q}$ cube-free, $A$ not a power of $\charfld(\finfldq{q})$ and $K = \ratsfuncqextn{q}{\varepsilon}$, where $\varepsilon$ is a root of $f(x)$.  Let $g$ be the genus of $K$.  Then $\deg(A) - 1 \leq g$.
\label{ABoundgivenG}
\end{Lem}

\begin{proof}
The genus of $K$ is $g=2\deg(A) - \deg(I) - 2$ by Proposition \ref{ListOfInvar}.  Rearranging this equation yields $\deg(I) = 2\deg(A) - g- 2 \leq \deg(A) - 1$, with the last inequality following from the remarks preceding the statement of Lemma \ref{ABoundgivenG}.  Rearranging this last inequality yields the claimed lower bound $\deg(A)-1$ for $g$, as desired.
\end{proof}

Another useful result is the determination of how finite places split (the $P$-signature) in the Shanks family.  This is well-known in the cyclic case, and is detailed below.  It is also a\ straightforward consequence of Theorem 8.2 of \cite{LRWWS}.  We specialize immediately to the Shanks family in Proposition \ref{FinPrimeSplit}, as this is useful for the computation of class numbers, given in Section \ref{sec:compres}.

\begin{Prop}
Let $K/\ratsfuncq{q}$ be a simple cubic function field with minimal polynomial given in standard form.  That is $f(Y) = Y^3 - \mathcal{A}Y + \mathcal{B}$ with $\mathcal{A},\mathcal{B}$ as in the proof of Proposition \ref{ListOfInvar}, is in standard form and $D = D(f) = (A^2 + 3A + 9)^2$.   Let $P$ be a polynomial in $\polyringq{q}$.  Then $K/\ratsfuncq{q}$ has $P$-signature
\begin{itemize}
\item $(1,1;1,1;1,1)$ if 
\begin{itemize}
\item $v_P(\mathcal{A}) = 0 < v_P(\mathcal{B})$, and $\mathcal{A}$ is a square modulo $P$;
\item $v_P(\mathcal{A}) = v_P(\mathcal{B}) = 0$, $v_P(D)$ is even, and the congruence $t^3 - \mathcal{A}t + \mathcal{B} \equiv 0 \pmod{P}$ has three solutions in \polyringq{q}.
\end{itemize}
\item $(1,3)$ if $v_P(\mathcal{A}) = v_P(\mathcal{B}) = 0$, and $v_P(D)$ is even, and the congruence $t^3 - \mathcal{A}t + \mathcal{B} \equiv 0 \pmod{P}$ has no solutions in \polyringq{q}.
\item $(3,1)$ if $P \mid \mathcal{A}$ and $P \mid \mathcal{B}$,
\end{itemize}
and no other cases occur.
\label{FinPrimeSplit}
\end{Prop}

The above invariants are easily computed for the non-Galois simple function fields by appealing to Theorem 8.2 of \cite{LRWWS}.   The statements for the Galois case have been included because of their relative elegance.   

\section{Fundamental Units and Regulator}\label{sec:Units}

In this section, we determine the fundamental units and regulator of the simple cubic function fields.
Determining the fundamental units of a unit rank two cubic function field can generally be done with Voronoi's Algorithm \cite{LeeSchYar,Sch1} or methods that rely on the infrastructure of the (ideal or divisor) class group \cite{Sch2}.  The units for Shanks simple cubic number fields arise quite naturally without having to use any time-costly algorithms; one fundamental unit is a root $\varepsilon$ of the defining polynomial and the other fundamental unit is $\varepsilon +1$.  
Shanks verified that these units were indeed fundamental by applying Godwin's criterion \cite{Godwin1,Godwin2}.  The units in the function field setting arise in precisely the same way, so we require a method to verify that the units  $\{\varepsilon, \varepsilon +k\}$ are indeed fundamental for both the Galois and non-Galois families of simple cubic function fields.

For this section, we now focus our attention back to arbitrary simple cubic fields arising from the polynomial
\begin{equation}\label{minpoly} x^3 - Ax^2 - (A + b)x - c \end{equation}
where $b \in \mathbb{F}_q$ and $c \in \mathbb{F}_q^*$.  By the construction in Section \ref{sec:FieldProps}, we know that these fields have the property that a root is a non-torsion unit and another independent unit can be obtained by adding a unit in $\mathbb{F}_q^*$.  We will now prove that when the function field is a simple cubic function field that the two units are fundamental units.  

While Godwin's criterion \cite{Godwin1} and its subsequent refinements by Cusick \cite{Cusick1,Cusick2} can be used either for determining whether units are fundamental or used constructively \cite{CusickComp} in number fields, we take a different approach to proving that the units  $\{\varepsilon, \varepsilon +k\}$ are fundamental for the function field setting.    First, we require some lemmas regarding the general case where the units under consideration are $\varepsilon$ and $\varepsilon + k$ in an arbitrary simple cubic function field.


\begin{Lem}
Let $K$ be a simple cubic function field with minimal polynomial (\ref{minpoly}).  Then there exists a root $\varepsilon$ of (\ref{minpoly}) that has positive degree.  Furthermore $\deg(\varepsilon) = \deg(A)$.
\end{Lem}

\begin{proof}
Let $\varepsilon, \varepsilon',\varepsilon''$ denote the roots of $x^3 - Ax^2 - (A + b)x - c$.  Since $\norm{\varepsilon} \in \finfldq{q}^*$, the degrees of the roots must add up to zero.  This rules out all the roots being negative.  Given the constraint on the norm, the only other possibilities are that each root has degree zero, or there exists a root with positive degree.  We now rule out the former.\\

Suppose all roots have degree zero.  We note that for any root $\tau \in \{\varepsilon, \varepsilon',\varepsilon''\}$,
\begin{equation}\tau^3 - A\tau^2 -(A+b)\tau - c = 0. \label{UnitRelation} \end{equation}  
Let $A = a_nt^n + a_{n-1}t^{n-1} + ... + a_0$, with $a_n \neq 0$.  For a fixed $\tau \in \{\varepsilon, \varepsilon',\varepsilon''\}$, we have $\tau \in \finfldq{q}\langle t^{-1} \rangle$, hence we can write $\tau = e_0t^0 + e_{-1}t^{-1} + \cdots$ with $e_0 \neq 0$.  Then $\deg(A) \geq \deg(\tau)$ and the leading term on the left-hand side of (\ref{UnitRelation}) is $a_ne_0^2 t^{n}$ arising from the $A\tau^2$ expansion.  But the coefficient of this leading term must be zero, hence $e_0^2=0$ or $a_n=0$.  In either case, this contradicts the assumption that $a_n \neq 0$ and $e_0 \neq 0$.  Hence there indeed exists a root of (\ref{minpoly}) with positive degree.  We label this root $\varepsilon$.

We now prove the last assertion, namely that $\deg(\varepsilon) = \deg(A)$.  First, suppose $m>n$.  Then the highest degree term in equation (\ref{UnitRelation}) is $t^{3m}$ arising from the $\varepsilon^3$ expansion.  But since the right-hand side of equation (\ref{UnitRelation}) is zero, then $e_m^3 = 0$, so $e_m = 0$, a contradiction.

Now, suppose that $m<n$.  Then repeating the argument we used to show that not all roots of (\ref{minpoly}) have degree zero yields the result on replacing $e_0$ with $e_m$.
Hence, it follows that $\deg(\varepsilon) = \deg(A)$, i.e.\ $m=n$.
\end{proof}


To simplify the subsequent discussion, always assume that a simple cubic function field is obtained by adjoining a root of (\ref{minpoly}) of positive degree to $\ratsfuncq{q}$, as adjoining any other root yields a conjugate field.
We now treat the Galois and non-Galois cases of simple cubic function fields separately for the next two propositions, even though the end result is identical in both cases.  First, suppose $\varepsilon$ is a root of $f(x) = x^3 - Ax^2 - (A + b)x - c$, where $f$ defines a non-Galois simple cubic function field. It remains to find expressions for the remaining roots in terms of $\varepsilon$.  We use the relation $ x^3 - Ax^2 - (A + b)x - c = (x-\varepsilon)(x-\varepsilon')(x-\varepsilon'')$.  It follows that $\varepsilon\varepsilon'\varepsilon'' = c$  and $\varepsilon'' = A - \varepsilon - \varepsilon'$.  Combining these two expressions we obtain $\varepsilon'(A - \varepsilon - \varepsilon') = c/\varepsilon$.  Rewriting this expression yields $\varepsilon'^2 + (A-\varepsilon)\varepsilon' + c/\varepsilon = 0$, a quadratic equation in $\varepsilon'$.  It has solutions 
\[\varepsilon' = \frac{-(A-\varepsilon) \pm\sqrt{(A-\varepsilon)^2 - 4c/\varepsilon}}{2}.\]
We label the solution with the $+$ sign before the square root as $\varepsilon'$, and the other solution as $\varepsilon''$.  These are the roots of the original polynomial in terms of $\varepsilon$.
We note that the degree of $\varepsilon'$ is $-\deg(\varepsilon) = -\deg(A)$, and the degree of $\varepsilon''$ is zero.  Note that in the case of a simple cubic function field, $\varepsilon$ is a unit with $\varepsilon \notin \finfldq{q}$, and that the same holds for $\varepsilon + k$.  Since these two units have the same degree and are not $\finfldq{q}^*$ multiples of each other, they must be independent.   The same observation is true for Shanks simple cubic function fields and the units $\varepsilon$ and $\varepsilon + 1$.

\begin{Prop}
Let $K$ be a non-Galois simple cubic function field defined by $ x^3 - Ax^2 - (A + b)x - c$. Then the regulator $R$ of $K$ is bounded above by $\deg(A)^2$.
\label{UpBdForRNonGal}
\end{Prop}

\begin{proof}

Since $\{\varepsilon, \varepsilon + k\}$ are two independent units that generate a subgroup of finite index of $\algints$, the regulator $R$ may be bounded by
\begin{eqnarray*}
R & \leq &\left| \det \left( \begin{array}{cc} \deg(\varepsilon) & \deg(\varepsilon') \\ \deg(\varepsilon + k) & \deg(\varepsilon +k)' \end{array} \right) \right|\\ 
    & = & \left| \det \left( \begin{array}{cc} \deg(A) & \deg(\varepsilon') \\ \deg(A) & \deg(\varepsilon' +k) \end{array} \right) \right| \\
    & = &\left| \det \left( \begin{array}{cc} \deg(A) & \deg\left(\frac{-(A-\varepsilon) + \sqrt{(A-\varepsilon)^2 - 4c/\varepsilon}}{2}\right) \\ \deg(A) & \deg\left(\frac{-(A-\varepsilon) +\sqrt{(A-\varepsilon)^2 - 4c/\varepsilon}}{2}+k\right) \end{array} \right) \right| \\
  & = & \left| \det \left( \begin{array}{cc} \deg(A) & - \deg(A) \\ \deg(A) &0 \end{array} \right) \right| \\
& = & \deg(A)^2,
\end{eqnarray*}
as desired.
\end{proof}

\begin{Prop}
Let $K$ be a Shanks simple cubic function field. Then the regulator $R$ of $K$ is bounded above by $\deg(A)^2$.
\label{UpBdForR}
\end{Prop}

\begin{proof}
The proof is virtually identical to the proof of Proposition \ref{UpBdForRNonGal}, except in Shanks fields, we always have $k=1$ and the roots of the Shanks polynomial are $\varepsilon$, $\varepsilon' = -1/(\varepsilon + 1)$ and $\varepsilon'' = -1/(\varepsilon' +1)$.

\end{proof}

\noindent Cusick \cite{Cusick3} proved that 
\[ R \geq \frac{1}{16}\log^2{(D/4)} \]
holds true for totally real cubic number fields with discriminant $D$.  Adapting this result to cubic function fields of unit rank two, we get the following result.

\begin{theorem}
Let $D$ be the discriminant of a cubic function field $K$ of unit rank $2$.  Then the regulator $R$ of $K$ satisfies
\[ R \geq \frac{1}{16}\deg^2{D} .\]
\label{LrBdForR}
\end{theorem}

\begin{proof}
The result follows in a similar way as the proof given in Cusick \cite{Cusick3}.  We simply note that our result is obtained by replacing $\log$ with $\deg$, and noting that the $\log 4$ factor becomes $\deg (4) = 0$ in the function field case.
\end{proof}

Note that when the parameters $A, c$, and $b$ are chosen so that (\ref{minpoly}) is nonsingular then $D(\varepsilon) = D(K)$; thus, it follows that
 \[ R \geq \frac{1}{16}\deg^2{(D(K))} = \frac{1}{16}(4\deg{(A)})^2 = \deg^2{(A)},\]
yielding $R = \deg(A)^2$ when combined with Proposition \ref{UpBdForRNonGal}.  When the curve is chosen to be singular, we get a small interval in which the regulator may lie. This does not happen frequently; ``randomly" chosen parameters $A$, $c$, and $b$ are likely to result in a nonsingular curve.  One way to force singularities to occur is by choosing $A$ to be a $p$-th power when $\textrm{char}(\mathbb{F}_q) = p$.  However, as previously mentioned, in this case an isomorphic function field may be considered by taking the $p$-th root of $A$ as a parameter for the minimal polynomial of $K$ instead.  

In fact, we can prove that for any simple cubic function field satisfying $\deg(I) \leq \deg(A) -1$ that $R=\deg(A)^2$, and there is no need to worry about a large gap between the upper and lower bound for $R$.  We give the general result first, followed by a corollary in the case of Shanks' family.

\begin{theorem}
 Let $K$ be a non-Galois simple cubic function field with discriminant $D(K)$, with minimal polynomial $f(x) = x^3 - Ax^2 -  (A+b)x - c$ and $\deg(I) \leq \deg(A) -1$.  Then $R = \deg(A)^2$ and $\{\varepsilon, \varepsilon+k\}$ are fundamental units.
 \label{FundUnitsProofGen}
\end{theorem}

\begin{proof}
We know from Theorems \ref{UpBdForR} and \ref{LrBdForR} that 
\[\deg(A)^2 \geq R \geq \frac{1}{16} \deg(D)^2. \]
Rewriting this inequality completely in terms of $A$ and $I$, we get

\begin{eqnarray*}
\deg(A)^2 \geq & R &  \geq  \frac{1}{16} (4\deg(A) - 2 \deg(I))^2 \\
\deg(A)^2 \geq & R &  \geq  \frac{1}{16} (16 \deg(A)^2 -16 \deg(A)\deg(I) + 4 \deg(I)^2) \\
\deg(A)^2 \geq & R &  \geq  \deg(A)^2 - \deg(A)\deg(I) + \frac{1}{4} \deg(I)^2
\end{eqnarray*}
Put $a = \deg(A)$ and $x= \deg(I)$.  Our inequality can now be rewritten, as $a^2 \geq R \geq a^2 - ax + 1/4 (x^2)$.  \\

\noindent Since $R$ must divide $ \deg(A)^2$, it is sufficient to show that $R > a^2/2$ in order to obtain the result.  The difference between the upper and lower bound for $R$ is $ax - \frac{1}{4} x^2$, so suppose for contradiction that the difference between the upper and lower bounds is at least $a^2/2$.

\noindent This yields
\begin{eqnarray*}
ax - \frac{1}{4} x^2 & \geq & \frac{a^2}{2}\\
4ax - x^2 & \geq & 2a^2 \\
x^2 - 4ax + 2a^2 & \leq & 0
\end{eqnarray*}
We denote the roots of the above quadratic by $\theta_1 = 2a + a/\sqrt{2}$ and $\theta_2 = 2a - a/\sqrt{2}$.  In order for the inequality $x^2 - 4ax + 2a^2  \leq  0$ to hold, we require either i):  $x-\theta_1 \leq 0$ and $x - \theta_2 \geq 0$ or ii):  $x-\theta_1 \geq 0$ and $x - \theta_2 \leq 0$.\\

\noindent Case i)  If $x - \theta_2 \geq 0$, then $x - 2a - a/\sqrt{2} \geq 0$.  Solving this inequality for $x$, we see that $x \geq a(2 - \frac{1}{\sqrt{2}}) = 1.29289 a$.  This contradicts the assumption that $x = \deg(I) \leq a -1$. \\

\noindent Case ii)  If $x-\theta_1 \geq 0$, then $x \geq 2a + \frac{a}{\sqrt{2}}$.  This also contradicts the assumption that $x=\deg(I) \leq a -1$.\\

\noindent In both cases, we have a contradiction.  Thus, it follows that $a^2 \geq R > a^2/2$.  The only number in the interval $[a^2,a^2/2)$ that divides $a^2$ is $a^2$ itself.  This completes the proof of the theorem.
\end{proof}

\begin{Cor}
Let $K$ be a Shanks simple cubic function field with discriminant $D(K)$, with minimal polynomial $f(x) = x^3 - Ax^2 -  (A+3)x - 1$ and $A^2 + 3A + 9$ cube-free.  Then $R = \deg(A)^2$ and $\{\varepsilon, \varepsilon+1\}$ are fundamental units.
\label{ShanksCor}
\end{Cor}

\begin{proof}
The proof proceeds exactly as in Theorem \ref{FundUnitsProofGen}, with the $A^2 + 3A + 9$ cube-free assumption yielding the bound $\deg(I) \leq \deg(A)-1$ as needed.
\end{proof}

The Shanks cases we have examined empirically with $A^2 + 3A + 9$ not cube-free may have singular points that do not invovle $p$-th powers.  Some of these fields, but not all of them, have non-trivial regulator less than $\deg(A)^2$.  The assumption that $A^2 + 3A + 9$ cube-free allows us to bound the index $I$ in terms of $A$; if $A^2 + 3A + 9$ is not cube-free, high values of $\deg(I)$ can and do occur, resulting in smaller regulators.  We discuss this and give examples in the next section.

\section{Computational Results:  Ideal Class Numbers}\label{sec:compres}

In this section, we proceed in an analogous fashion to Shanks' original work \cite{ShanksSimp}.  Namely, we give the results of our ideal class number computations of these fields.  We divided this section into two parts.  The first part covers some general results of our computations, a brief recall of the Hasse-Weil bounds and other miscellaneous results.  We detail the solution to finding all Galois simple cubic function fields with ideal class number $1$ in the second part.   We used Magma \cite{Magma} and C/C++ with NTL \cite{Shoup1} for our computations.  Computations were carried out on a MacBook Pro with a 2.16 GHz Intel Core 2 Duo processor and 2GB of RAM, and on a multi-processor machine with four 2.8 GHz Pentium 4 processors running Linux with 4 GB of RAM.

\subsection{Results for Shanks families over $\finfldq{5}$ and $\finfldq{7}$}

One of our main tools in this section and subsequent sections is the Hasse-Weil bound (\cite{Sticht}, Theorem 5.2.3), an equivalent version of which we give below.

\begin{theorem}[Hasse-Weil]
Let $K/\ratsfuncq{q}$ be a function field over $\ratsfuncq{q}$ of genus $g$.  Then the divisor class number $h$ satisfies
\begin{equation*}(\sqrt{q}-1)^{2g} \leq h \leq (\sqrt{q}+1)^{2g} .\end{equation*}
\label{HW}
\end{theorem}
We also give a structure theorem on the ideal class group of $\algints$.  The proof is precisely the same as the one given in Shanks \cite{ShanksSimp}.

\begin{theorem}
Let $K/\ratsfuncq{q}$ be a Galois simple cubic function field over $\ratsfuncq{q}$ of genus $g$.  Then the ideal class number $h'$ is of the form $a^2 + 3b^2$ for some pair of integers $(a,b)$.
\label{ShanksGrpStruct}
\end{theorem}

We note that Theorem \ref{ShanksGrpStruct} also holds when $h'$ is replaced by $h$ when $A^2+3A+9$ is cube-free, since then the two quantities differ by a square, namely $R = \deg(A)^2$. 

Our strategy for computing class numbers was as follows.  Ideally, one would use the procedure given in \cite{SchStein1}, namely use the truncated Euler product of the zeta function of $K$ to compute an approximation $E$ of $h'$ and an integer $L$ such that $|h-E|< L^2$, and then use Shanks' baby step giant step (BSGS) to find the correct value of $h'$ in the interval $(E-L^2,E+L^2)$.  Unfortunately, efficient ideal arithmetic and BSGS methods for cubic function function fields of unit rank two have not been developed, other than the case of purely cubic fields \cite{FontLandSch}.  While such methods would apply for Shanks' family when $q\equiv 1 \pmod{3}$, as such fields have a purely cubic model (see \cite{LeeSchYar}, Lemma 2.1), we would still be stuck in the case when $q \equiv -1 \pmod{3}$.

Our approach was to experiment, and make exclusive use of the truncated Euler product when computing these class numbers.  While the Euler product does converge fairly slowly, as in the case of simple cubic number fields, we only made use of the Euler product for very small values of $\deg(A)$.  We experimented with varying the size of the product depending on the size of $\deg(A)$.  We adopted Shanks \cite{ShanksSimp} remark ``a great deal of accuracy is not needed here since $h$ is an integer" as our philosophy, especially in the case of small values of $\deg(A)$ when looking for class number $1$ fields in Subsection \ref{ssec:clasno1}.  One advantage of using the ``Euler product only" approach is its simplicity.

We briefly describe the method of using the Euler product to compute class numbers.  A more thorough description can be found in \cite{SchStein1} and \cite{SchStein2}.  The $\zeta$-function of a function field $K$ of genus $g$ is defined by 
\[\zeta(s,K) = \sum_{\mathfrak{A}}  \frac{1}{N(\mathfrak{A})^s} \;\;\;\;\;\;\;  \mathrm{Re}(s) > 1,\]
where the summation is over all integral divisors $\mathfrak{A}$ of $K$ and $\mathrm{Re}(s)$ denotes
the real part of the complex variable $s$.  It is standard in the function field setting to write $u=q^{-s}$ and write $\zeta(s,K) = Z(u,K)$.  The function $Z(u,K)$ has an Euler product formula, and is in fact a rational function in $u$.  Furthermore, $Z(u,K)$ can be written as $Z(u,K) = Z(u,k(X))L(u,K)$ where $L$ is the $L$-polynomial
\[L(u,K) = \prod_{i=1}^{2g}(1-\alpha_iu).\]
The $L$-polynomial satisfies a functional equation, and the value of $L(1,K)$ is the divisor class number $h$.  In fact, one can show that 
\begin{equation}h = L(1,K) = \frac{q^{g+2}}{(q-x_1)(q-x_2)} \prod_{P}f(P,1/q), \label{hEulProd}\end{equation}
where the product ranges over all monic irreducible polynomials $P$, and $(x_1,x_2)$ and the quantities $f(P,1/q)$ are determined by how the infinite place of $K$ splits, and how the finite places $P$ split in $K$, respectively.  The latter part can be explicitly determined using Proposition \ref{FinPrimeSplit}.  Approximating $h$ is then a matter of deciding where to truncate the Euler product in Equation (\ref{hEulProd}).  Once this $h$ approximation is found, $h'$ is found via the relation $h = Rh'$. 

We present a sample of some of the results of our computations for index-free Galois simple cubic fields in Table \ref{ShanksFldsNoIndex}.  Examples with non-trivial index were also computed and are omitted for brevity.  Column $1$ denotes the order of the underlying finite field.  Column 2 gives the value of the parameter $A$ determining the minimal polynomial of the Shanks simple cubic function field.  
The column denoted $h'$ gives the ideal class number of the field.  The columns denoted \# of Split $P$ and \# of Inert $P$ denote the number of split primes and the number of inert primes up to the truncation point, respectively.  The column labelled ``Max.\ Measure" denotes how close we get to upper end of the Hasse-Weil interval for each ideal class number, and is simply the value of 
\[\frac{h'-\left\lceil \frac{(\sqrt{q}-1)^{2g}}{R}\right\rceil}{\left\lfloor \frac{
(\sqrt{q}+1)^{2g}}{R}\right\rfloor - \left\lceil \frac{(\sqrt{q}-1)^{2g}}{R} \right\rceil}, \]
where $R = \deg(A)^2$.  

The truncation point for $\deg(A) = 2,3$ was 40000 and $\deg(A) = 4$ was 100000, over $\finfldq{5}$.    The truncation point for $\deg(A) = 2$ was 40000 and for $\deg(A) = 3$ was 100000, over $\finfldq{7}$.    
The computations in these tables we done using C/C++ with NTL \cite{Shoup1}, and were verified with Magma \cite{Magma}.

\begin{center}
\begin{longtable}{|c|c|c|c|c|c|} 
\caption{Shanks simple cubic function fields over $\ratsfuncq{q}$ with trivial index $I$}\label{ShanksFldsNoIndex} \\

\hline \multicolumn{1}{|c|}{\,\,$q$\,\,} & \multicolumn{1}{|c|}{$A$} & \multicolumn{1}{|c|}{$h'$} & \multicolumn{1}{|c|}{\# of Split $P$} & \multicolumn{1}{|c|}{\# of Inert $P$} & \multicolumn{1}{|c|}{Max.\ Measure} 
\\ \hline \endfirsthead

\multicolumn{6}{c}%
{{\bfseries \tablename\ \thetable{} -- continued from previous page}} \\
\hline \multicolumn{1}{|c|}{\,\,$q$\,\,} & \multicolumn{1}{|c|}{$A$} & \multicolumn{1}{|c|}{$h'$} & \multicolumn{1}{|c|}{\# of Split $P$} & \multicolumn{1}{|c|}{\# of Inert $P$} & \multicolumn{1}{|c|}{Max.\ Measure} 
\\ \hline
\endhead

\hline \multicolumn{6}{|c|}{{Continued on next page}} \\ \hline
\endfoot

\hline 
\endlastfoot
\,\,$5$\,\,  & $t^2$ & $12$ & 13419 & 26579 & 0.42308 \\ \cline{2-6}
& $t^2+1$& $4$ &13324 & 26675 & 0.11538 \\ \cline{2-6} 
& $t^2+2$& $12$ &13387 & 26611 & 0.42308 \\ \cline{2-6}
& $t^2+3$& $4$ &13329 & 26670 & 0.11538 \\ \cline{2-6}
& $2t^2$& $3$ &13364 & 26634 & 0.07692 \\ \cline{2-6}
& $2t^2 + 1$& $4$ &13190 & 26809 & 0.11538 \\ \cline{2-6}
& $2t^2 + 2$& $3$ &13388 & 26610 & 0.07692 \\ \cline{2-6}
& $2t^2 + 3$& $13$ &13296 & 26703 & 0.5000 \\ \cline{2-6}
& $t^3$ & $73$ & 13377 & 26622 & 0.05464 \\ \cline{2-6}
& $t^3 + 3$ & $52$ & 13258 & 26741 & 0.03892 \\ \cline{2-6}
& $t^3 + t$ & $39$ & 13393 & 26605 & 0.02919 \\ \cline{2-6}
& $t^3 + t+1$ & $256$ & 13387 & 26612 & 0.19162 \\ \cline{2-6}
& $t^3 + t + 3$ & $63$ & 13275 & 26723 & 0.04716 \\ \cline{2-6}
& $t^3 + 2t$ & $27$ & 13355 & 26643 & 0.02021 \\ \cline{2-6}
& $t^3 + 2t+1$ & $81$ & 13258 & 26739 & 0.06063 \\ \cline{2-6}
& $t^3 + 2t+3$ & $111$ & 13405 & 26593 & 0.08308 \\ \cline{2-6}
& $t^4$ & $8112$ & 33608 & 66390 & 0.09841 \\ \cline{2-6}
& $t^4+1$ & $592$ & 33266 & 66733 & $7.182 \times 10^{-3}$ \\ \cline{2-6}
& $t^4+2$ & $768$ & 33161  & 66837 & $9.317 \times 10^{-3}$ \\ \cline{2-6}
& $t^4+3$ & $976$ & 33183 & 66816 & $0.0118$ \\ \cline{2-6}
& $t^4+4$ & $832$ & 33559 & 66440 & $0.0101$ \\ \cline{2-6}
& $t^4 + t$ & $468$ & 33390 & 66608 & $5.677 \times 10^{-3}$ \\ \cline{2-6}
& $t^4+t+1$ & $1812$ &  33383& 66615 & 0.02198 \\ \cline{2-6}
& $t^4+t+3$ & $868$ &  33232& 66767 & 0.01053 \\ \cline{2-6}
& $t^4 + t^2$ & $336$ & 33143 & 66855 & $4.076 \times 10^{-3}$ \\ \cline{2-6}
\hline
\,\,$7$\,\, &$t^2$ & 13  & 13350  & 26648 & 0.2619\\ \cline{2-6}
& $t^2+2$& 12 & 13362 & 26635 & 0.2381\\ \cline{2-6}
& $t^2+4$& 9 & 13307 & 26689 & 0.2143\\ \cline{2-6}
& $3t^2+3t+3$& 13 & 13365 & 26633 & 0.2619\\ \cline{2-6}
& $3t^2+3t+6$& 9 & 13306 & 26690 & 0.2143\\ \cline{2-6}
& $3t^2+4t+1$& 12 & 13369 & 26628 & 0.2381\\ \cline{2-6}
& $3t^2+4t+3$& 13 & 13364 & 26634 & 0.2619\\ \cline{2-6}
& $3t^2+5t$& 12 & 13307 & 26689 & 0.2381\\ \cline{2-6}
& $t^3+1$& 441 & 32949 & 67049 & 0.1257\\ \cline{2-6}
& $t^3+2$& 144 & 33433 & 66565 & 0.0399\\ \cline{2-6}
& $t^3+t+2$& 324 & 33238 & 66758 & 0.0922\\ \cline{2-6}
& $t^3+t+5$& 225 & 33126 & 66871 & 0.0633\\ \cline{2-6}
& $t^3+3t$& 117 & 33313 & 66684 & 0.0321\\ \cline{2-6}
& $t^3+3t + 2$& 729 & 33316 & 66680 & 0.2089\\ \cline{2-6}
& $t^3+5t + 1$& 252 & 33316 & 66680 & 0.0711\\ \cline{2-6}
& $2t^3+2t^2+6t + 1$& 576 & 33308 & 66690 & 0.1647\\ \cline{2-6}
\end{longtable}
\end{center}






Following \cite{SchStein1}, we can estimate the error of the approximations we computed with the parameters specified below.  We write $h$ in ``truncated Euler product form":
\begin{equation}
h = E' \cdot e^{B},  \label{truncEulProdForm}
\end{equation}
where $E'$ and $B$ are real numbers.  Note that $B = \log h - \log(E')$.  In order to get good approximations, we need a sharp upper bound $\psi \in \mathbb{R}$ on $|B|$.  Furthermore, if $\psi$ is noticeably smaller than $1$, then $|e^B-1| < e^\psi -1$ and we put $E = \mathrm{round}(E')$ and $L = \left\lceil \sqrt{E'(e^\psi -1)+\frac{1}{2}}\right\rceil$.  From this, it follows that $|h-E| \leq L^2$ as desired.

The estimate we use for our error is the ``second estimate" given in \cite{SchStein1}, which is specified by
\[E(\lambda,K) =  \frac{q^{g+2}}{(q-x_1)(q-x_2)} \prod_{P, \deg(P)=\nu \leq \lambda}f(P,1/q).\]
and
\[B(\lambda,K) = \log \,\,\,\,\,\,  \prod_{P, \deg(P)=\nu > \lambda} \frac{q^{2\nu}}{(q^\nu - z_1(P))(q^\nu - z_2(P))},  \]
where $z_1,z_2$ are determined by how each prime $P$ in the product splits in the field $K$.  An upper bound on $B_2(\lambda,K)$ is
\begin{eqnarray*}
\psi(\lambda,K) = \frac{2g}{\lambda+1}q^{-\frac{\lambda+1}{2}} & + & \frac{2g+4}{\lambda+2}\frac{\sqrt{q}}{\sqrt{q}-1}\frac{q}{q-1}q^{-\frac{\lambda+2}{2}} + \frac{2}{\lambda+1}q^{-(\lambda+1)} \\
& + & \frac{2}{\lambda+1}\frac{q}{q-1}q^{-(\lambda+1)} (q^{\frac{\lambda+1}{\epsilon_\lambda}}-1),
\end{eqnarray*}
where $\epsilon_\lambda = 2$ if $\lambda$ is odd and $\epsilon_\lambda = 3$ if $\lambda$ is even.  In practice, we truncated our products not by the degree value to be computed in the infinite product in equation (\ref{hEulProd}), but by the total number of polynomials in the product.  For instance, our computation of class numbers over $\finfldq{5}$ and $\deg(A)=2$ used the truncation point of $40000$ polynomials.  The last polynomial in the product has degree 8, but since we do not use all polynomials of degree $8$ in the product, we use a value of $\lambda=7$ to compute $\psi$ above.  Because of this, and because we know the value of the regulator in advance, our approximations for $h'$ as derived from the approximations of $h$ are more accurate than indicated by $\psi$, but this still gives us a good idea as to the accuracy of our approximations.  Upper bounds for $\psi$ for the various parameters $q$, $g$ and $\lambda$ that were used for computing the values in Table \ref{ShanksFldsNoIndex} are given in Table \ref{PsiTable}.  The last column in Table \ref{PsiTable} denotes the maximum value for $L$ found at this genus level for these fields, which is determined using the largest approximation found with the given parameters, and the formula for $L$ previously specified.

\begin{center}
\begin{longtable}{|c|c|c|c|c|} 
\caption{Values for $\psi$ and $L$ given $q$,$g$ and $\lambda$:  Fields with trivial index $I$}
\label{PsiTable} \\

\hline \multicolumn{1}{|c|}{\,\,$q$\,\,}  & \multicolumn{1}{|c|}{$g$} &  \multicolumn{1}{|c|}{$\lambda$}& \multicolumn{1}{|c|}{$\psi$} & \multicolumn{1}{|c|}{Maximum $L$}
 \\ \hline \endfirsthead

\multicolumn{5}{c}%
{{\bfseries \tablename\ \thetable{} -- continued from previous page}} \\
\hline \multicolumn{1}{|c|}{\,\,$q$\,\,}  & \multicolumn{1}{|c|}{$g$} &  \multicolumn{1}{|c|}{$\lambda$}& \multicolumn{1}{|c|}{$\psi$} & \multicolumn{1}{|c|}{Maximum $L$}
\\ \hline
\endhead

\hline \multicolumn{5}{|c|}{{Continued on next page}} \\ \hline
\endfoot

\hline 
\endlastfoot
\,\,$5$\,\, &  2 & 7  & $2.730717 \times 10^{-3}$ &  1 \\ \cline{2-5} 
&  4 & 8  & $1.522115 \times 10^{-3}$ &  3\\ \cline{2-5} 
&  6 & 8  & $2.129576 \times 10^{-3}$ & 17 \\ \hline
\,\,$7$\,\, &  2 & 6  & $1.448719 \times 10^{-3}$ & 1 \\ \cline{2-5} 
&  4 & 6  & $2.431095 \times 10^{-3}$ & 5 \\ \cline{2-5} 
\end{longtable}
\end{center}

Given these $\psi$ and $L$ values, one would expect very accurate approximations for our computations up to genus 6.  Our computations largely bear this out.  For instance, the chosen truncation point for $\finfldq{5}$ result in correct, exact values for $\deg(A)=2$.  In the $\deg(A)=3$ case, we get exact values plus a few instances where rounding to the next nearest integer results in the correct value for $h'$.  The truncation point for $\deg(A)=4$ produces additional error, but in the worst cases, the value of $h'$ is off by less than $5$.  This is largely expected due to the fact that larger and larger class approximations and class numbers are obtained as the genus increases, by equation (\ref{truncEulProdForm}).  Similar results were obtained for Shanks fields over $\finfldq{7}$.

The total running time of finding the approximation of $h$ is $O(q^{1/4})$ for genus $g=1,2$.  For $g \geq 3$, determining an approximation $E$ of $h$ using the methods in \cite{SchStein1} requires the computation of $O(q^\lambda)$ values $z_1(P),z_2(P)$ that determine $(P,1/q)$ in equation (\ref{hEulProd}) above for a given fixed $\lambda$.  For small genus and large $q$, $E = O(q^g)$.  While we work with smaller $q$ in this paper, this still provides a rough idea as to the complexity of the algorithm. Complete details can be found in \cite{SchStein1}.  
A table of our timings is given below.  We precomputed the irreducible polynomials and stored them in a table before computing the class number for each seed $A$ for a given $\deg(A)$ value.  As indicated above, we opted for extremely precise approximations which did impact the running time more adversely than in an approximation/BSGS approach.  This is still the subject of future research.

\begin{center}
\begin{longtable}{|c|c|c|c|c|} 
\caption{Timings for fields with trivial index $I$}
\label{TimingsTableNoIndex} \\

\hline \multicolumn{1}{|c|}{\,\,$q$\,\,}  & \multicolumn{1}{|c|}{$g$} & \multicolumn{1}{|c|}{\# of $A$ processed} &  \multicolumn{1}{|c|}{Total time} &  \multicolumn{1}{|c|}{Time per $A$}  \\ \hline \endfirsthead

\multicolumn{5}{c}%
{{\bfseries \tablename\ \thetable{} -- continued from previous page}} \\
\hline \multicolumn{1}{|c|}{\,\,$q$\,\,}  & \multicolumn{1}{|c|}{$g$} & \multicolumn{1}{|c|}{\# of $A$ processed}  &\multicolumn{1}{|c|}{Total time} &  \multicolumn{1}{|c|}{Time per $A$}\\ \hline
\endhead

\hline \multicolumn{5}{|c|}{{Continued on next page}} \\ \hline
\endfoot

\hline 
\endlastfoot
\,\,$5$\,\, &  2 & 100 & 1 hour, 23 min, 23 sec & 50 seconds \\ \cline{2-5} 
&  4 & 460  & 9 hours, 27 min, 36 sec & 74 seconds\\ \cline{2-5} 
&  6 & 2380  &3 days, 15 hours & 132 seconds\\ 
&    &   & 19 min, 34 seconds &   \\  \hline
\,\,$7$\,\, &  2 & 126 & 1 hour, 14 min, 14 sec & 35 seconds \\ \cline{2-5} 
&  4 & 1176 & 1 day, 11 hour& 109 seconds \\  
&    &   & 45 min, 30 seconds &   \\  \hline
\end{longtable}
\end{center}

We note that even though the regulator of the Galois simple cubic function fields are small, the ideal class numbers are also quite small.  As seen in Table~\ref{ShanksFldsNoIndex}, the ideal class numbers do not get anywhere close to the upper bound of the Hasse-Weil interval, so these function fields are nowhere near being maximal function fields.  Maximal function fields, i.e.\ fields where $h = (\sqrt{q}+1)^{2g}$, do exist for each $\finfldq{q^2}$.  One such example is the Hermitian function field, defined by $H = \finfldq{q^2}(t,y)$ with $y^q + y = t^{q+1}$ (See Stichtenoth \cite{Sticht}, Lemma 6.4.4).  The phenomenon of small class numbers also presents itself for simple cubic number fields, as noted by Shanks \cite{ShanksSimp}.  However, in the case of simple cubic number fields, Shanks notes that this small/moderate class number size is due to the fact that 2 and 3 are cubic non-residues for all $P=a^2+3a+9$ under his consideration.  We do not have a similar explanation in the function field case.

As previously mentioned, if the cube-free assumption on $A^2+3A+9$ is relaxed, the degree of the index $I$ is no longer bounded from above by $\deg(A) - 1$.  In these cases, larger indices can occur, resulting in regulators that are strictly smaller than $\deg(A)^2$.  The first example we were able to find was over $\finfldq{13}$, namely the Shanks field $K$ where $A = t^3+1$.  In this case, $\deg( D(K) ) = 6$ but $\deg(D(f)) = 12$ and  $D(f) = t^6 \cdot D(K)$.  The degree of the index is the same as the degree of $A$, and we obtain a regulator of $R=3$, instead of $R=9$.     A small sample of other examples found over $\finfldq{7}$ are given in Table \ref{HighIndexLowReg}.

\begin{center}
\begin{longtable}{|c|c|c|c|c|} 
\caption{A sample of Shanks simple cubic function fields over $\ratsfuncq{7}$ with large index $I$ and regulator $R < \deg(A)^2$}\label{HighIndexLowReg} \\

\hline \multicolumn{1}{|c|}{\,\,$q$\,\,}  & \multicolumn{1}{|c|}{$\deg(I)$} &  \multicolumn{1}{|c|}{$A$}& \multicolumn{1}{|c|}{$h'$} & \multicolumn{1}{|c|}{$R$} 
 \\ \hline \endfirsthead

\multicolumn{5}{c}%
{{\bfseries \tablename\ \thetable{} -- continued from previous page}} \\
\hline \multicolumn{1}{|c|}{\,\,$q$\,\,}  & \multicolumn{1}{|c|}{$\deg(I)$} & \multicolumn{1}{|c|}{$A$}&\multicolumn{1}{|c|}{$h'$} &\multicolumn{1}{|c|}{$R$}
\\ \hline
\endhead

\hline \multicolumn{5}{|c|}{{Continued on next page}} \\ \hline
\endfoot

\hline 
\endlastfoot
\,\,$7$\,\, &  3 & $t^3+5$  & 3 & 3  \\ \cline{3-5}
&  & $t^3+t^2 + 5t + 4$ & 3 & 3  \\ \cline{3-5}
&  & $t^3+2t^2 + 6t + 5$& 3 & 3  \\ \cline{3-5}
&  & $2t^3+ 5$ & 1 & 3  \\ \cline{3-5}
&  & $2t^3+ 6$ & 1 & 3  \\ \cline{3-5}
&  & $2t^3 +t^2 +6t+ 3$ & 1 & 3  \\ \cline{3-5}
\end{longtable}
\end{center}

The value of $A^2 + 3A +9$ being not cube-free does not guarantee a regulator value smaller than the upper bound $\deg(A)^2$.  For example, the $A$ values $t^3+3t+5$ and $t^3+3t+6$ yield $R=9$ even though $A^2 + 3A +9$ is not cube-free, so having a cube factor in $A^2 + 3A +9$ is not sufficient to guarantee a value for $R$ less than $\deg(A)^2$.

\subsection{Shanks Simple Cubic Fields with Class Number One}\label{ssec:clasno1}

As an application of our computational work and the Hasse-Weil bound, in this section we determine all Galois simple cubic function fields with class number~$1$ and $A^2+3A+9$ cube-free.  First, we note that from Lemma \ref{ABoundgivenG}, it follows that $\deg(A) - 1 \leq g \leq 2\deg(A) -2$.  It is not hard to show that the lower bound in Theorem \ref{HW}
\begin{equation}(\sqrt{q}-1)^{2(\deg(A)-1)}/\deg(A)^2    \label{HWLowBd}\end{equation}
is an increasing function in the variable $\deg(A)$.  Hence for a fixed $q$, there is some value $X$ of $\deg(A)$ such that equation (\ref{HWLowBd}) is greater than $1$ for all values of $a < X$.  This limits the potential candidates for $A$ yielding fields with class number~1, as seen in the following proposition, which is a straightforward application of Lemma \ref{ABoundgivenG} and Theorem \ref{HW}.

\begin{Prop}
Let $A \in \finfldq{q}$ with $A^2+3A+9$ cube-free, and let $K/\ratsfuncq{q}$ be a Shanks simple cubic function field defined by the minimal polynomial $x^3 - Ax^2 - (A+3)x - 1$.  Then if $K$ has ideal class number equal to $1$, then either $q=5$ and $\deg(A)\leq 13$ or $q=7$ and $\deg(A)\leq 3$ and there are no other possibilities.
\label{ClasNo1Candidates}
\end{Prop}

We can refine the results from Proposition \ref{ClasNo1Candidates} in a more effective way by using the genus formula, the Hasse-Weil bound and the fact that $\deg(I) \leq \deg(A)-1$.  This enables a more targeted search for fields with a specified index and $\deg(A)$ value.  The class number of these fields can then be computed using Magma \cite{Magma}.  We do this for $\finfldq{5}$ only, since $\deg(A)$ is small for class number~1 fields over $\finfldq{7}$.  We also exclude the trivial case where $\deg(A)=1$ from the proposition below.

\begin{Prop}
Let $A \in \finfldq{5}$ with $A^2+3A+9$ cube-free and $\deg(A)>1$, and let $K/\ratsfuncq{5}$ be a Shanks simple cubic function field defined by the minimal polynomial $x^3 - Ax^2 - (A+3)x - 1$.  Then if $K$ has ideal class number equal to $1$, then 
\begin{itemize}
\item If $\deg(I) = 0$, then $2 \leq \deg(A) \leq 4$.
\item If $\deg(I) = 1$, then $2 \leq \deg(A) \leq 5$.
\item If $\deg(I) = 2$, then $3 \leq \deg(A) \leq 6$.
\item If $\deg(I) = 3$, then $4 \leq \deg(A) \leq 7$.
\item If $\deg(I) = 4$, then $5 \leq \deg(A) \leq 7$.
\item If $\deg(I) = 5$, then $6 \leq \deg(A) \leq 8$.
\item If $\deg(I) = 6$, then $7 \leq \deg(A) \leq 9$.
\item If $\deg(I) = 7$, then $8 \leq \deg(A) \leq 9$.
\item If $\deg(I) = 8$, then $9 \leq \deg(A) \leq 10$.
\item If $\deg(I) = 9$, then $10 \leq \deg(A) \leq 11$.
\item If $\deg(I) = 10$, then $\deg(A) = 11$.
\item If $\deg(I) = 11$, then $\deg(A) = 12$.
\item If $\deg(I) = 12$, then $\deg(A) = 13$.
\end{itemize}
\label{ClasNo1CandidatesRefined}
\end{Prop}

From our computations and Propositions \ref{ClasNo1Candidates} and \ref{ClasNo1CandidatesRefined}, we obtain the following theorem.  We include the case where $\deg(A)=1$ for any $q$, yielding genus $0$ fields, for the sake of completeness.\\

\begin{theorem}
Let $A \in \finfldq{q}$ with $A^2+3A+9$ cube-free, and let $K/\ratsfuncq{q}$ be a Shanks simple cubic function field defined by the minimal polynomial $x^3 - Ax^2 - (A+3)x - 1$.  Then $K$ has ideal class number equal to $1$ if and only if either $\deg(A)=1$ or $A$ is one of the following possibilities given in Table \ref{ClasNo1table}. 
\end{theorem}
\begin{center}
\begin{longtable}{|c|c|c|} 
\caption{Shanks simple cubic function fields with class number 1}\label{ClasNo1table} \\

\hline \multicolumn{1}{|c|}{\,\,$q$\,\,} & \multicolumn{1}{|c|}{$\deg(I)$}  & \multicolumn{1}{|c|}{$A$} 
 \\ \hline \endfirsthead

\multicolumn{3}{c}%
{{\bfseries \tablename\ \thetable{} -- continued from previous page}} \\
\hline \multicolumn{1}{|c|}{\,\,$q$\,\,} & \multicolumn{1}{|c|}{$\deg(I)$} & \multicolumn{1}{|c|}{$A$} 
\\ \hline
\endhead

\hline \multicolumn{3}{|c|}{{Continued on next page}} \\ \hline
\endfoot

\hline 
\endlastfoot
\,\,$7$\,\, & 1 & $t^2+6$   \\ \cline{3-3}
& & $t^2 + t + 1$   \\ \cline{3-3}
& & $t^2 + 2t$    \\ \cline{3-3}
& & $t^2 + 3t + 3$     \\ \cline{3-3}
& & $t^2 + 4t +3$     \\ \cline{3-3}
& & $t^2 + 5t$  \\ \cline{3-3}
& & $t^2 + 6t +1$     \\ \cline{3-3}
& & $2t^2 +6$    \\ \cline{3-3}
& & $2t^2 +t$    \\ \cline{3-3}
& & $2t^2 +2t+3$    \\ \cline{3-3}
& & $2t^2 +3t+1$    \\ \cline{3-3}
& & $2t^2 +4t+1$    \\ \cline{3-3}
& & $2t^2 +5t+3$    \\ \cline{3-3}
& & $2t^2 +6t$    \\ \cline{3-3}
& & $3t^2 +5$    \\ \cline{3-3}
& & $3t^2 +t + 1$    \\ \cline{3-3}
& & $3t^2 +2t +3$    \\ \cline{3-3}
& & $3t^2 +3t +4$    \\ \cline{3-3}
& & $3t^2 +4t +4$    \\ \cline{3-3}
& & $3t^2 +5t +3$    \\ \cline{3-3}
& & $3t^2 +6t +1$    \\ \cline{3-3}
\end{longtable}
\end{center}
\label{ClasNo1ListFinal}

For smaller values of $\deg(A)$, we used the truncated Euler product to compute class numbers of the Shanks simple cubic function fields.  For larger values of $\deg(A)$, we identified class number $1$ fields in the following way.  We used a separate C++/NTL program to identify higher index fields for large values of $\deg(A)$ for $\finfldq{5}$, and then computed the class number of these candidates using Magma to find the class number 1 fields.  

As an experiment, we tested our C++ implementation on some of the $\finfldq{5}$ candidates indicated in Proposition \ref{ClasNo1CandidatesRefined}.  For instance, in searching for potential class number~1 candidates amongst the $\deg(I) = 10, \deg(A)=11$ fields, we found all such Shanks fields had $h'=76176$.  Our implementation, using a truncation point of $140000$ ($\lambda = 8$ in this case), gets within $25$ of the actual value for each of these fields.  In some cases, we are off by as little as $2$ with this same truncation point.  This is better than predicted by the $\psi$ value in these cases:  with the above parameters and $\lambda = 8$, we get an $L$ value of $176$, so $|h-E| < 30976$.  Again, since we are interested in the ideal class number $h'$, dividing the previous inequality by $R=121$ yields $|h' - E^*| < 256$, where $E^* = E/R$.  Since our product approximations are truncated at the $n$-th polynomial, and not after all polynomials of degree $8$, we do much better than the error bounds suggest.

Unlike Shanks \cite{ShanksSimp} and Lettl \cite{Lettl}, we did not restrict our computations to prime values of $A^2+3A+9$.  Rather, our cube-free assumption is stronger and yields some extra fields with class number 1.  We also found some class number~1 Shanks function fields with the cube-free restriction relaxed (see Table \ref{HighIndexLowReg}), but enumerating these fields with the cube-free assumption lifted is an open problem.



\section{Open Problems and Future Work}\label{sec:concl}

A number of open problems arise from this work.   Niklasch and Smart  \cite{quartic_exc} extended Nagell's study of exceptional units to quartic fields.   Since the infinite place can have arbitrary splitting type, many more extension degrees can permit exceptional units.  Thus, function fields display a bit more flexibility than their number field counterparts in this respect.  One may also investigate quintic and sextic extensions; these have been studied via Gaussian periods by Lehmer \cite{Lehmer1} and Schoof and Washington \cite{SchoofWash}.  A generalization of these results to the function field setting is the subject of future investigation.  A key connection between elliptic curves and the $2$-part of the class group of simple cubic number fields was established by Washington \cite{Wash1}, and is also the subject of future research.

To prove that the set of units that arose from the roots of the defining polynomial were fundamental in the number field case, Shanks used Godwin's Criterion \cite{Godwin1}.  This result, and subsequent refinements by Cusick \cite{Cusick1,Cusick2,CusickComp}, do not seem to generalize easily to function fields.  Using Godwin's criterion involves finding the minimum value of a certain positive definite ternary quadratic form.  This is problematic due to the non-Archimedean nature of the absolute value in the function field setting, necessitating our approach.  Another way to generate such simple fields in the number field setting is to use modular curves, as Washington does in the quartic case \cite{Wash2}.  It is also noted in the aforementioned work that the simple cubic fields also arise in this manner.  Generalizing this approach to the function field setting using Drinfeld modular forms remains an open problem.

The computational aspects of this work may also be improved.  However, efficient ideal arithmetic is only available for purely cubic function fields \cite{bauer_cubic,FontLandSch}.  If $q \equiv 1 \pmod{3}$ then a Galois cubic function field has a purely cubic model.  Our computations over $\mathbb{F}_7$ could have been changed to a purely cubic model in order to use more efficient ideal arithmetic, but we chose to make explicit use of the Shanks model.  Since general divisor arithmetic is available \cite{hess, mak}, it would be nice to know how our computational results compare with using the more typical reduction methods to find the regulator \cite{Sch2}.   

Finally, we have omitted detailed computations for the case of non-Galois simple cubic function fields.  This will be the subject of a future paper, along with the development of criteria for determining whether two non-Galois families are isomorphic.  The isomorphism question for Galois cubics is discussed in Kersten and Michali\v{c}ek \cite{KerstMicha}.

\bibliographystyle{amsplain}
\bibliography{exunitsreffile}
\end{document}